\documentclass[11pt,bezier]{article}
\usepackage{amsmath,amssymb,amsfonts,euscript,graphicx}

\title{\LARGE\bf   On Commutative Rings Whose Prime Ideals   Are Direct Sums of Cyclics\thanks {The research
 of first author   was in part supported by
a grant from IPM (No. 90160034).}
\thanks
 {{\it Key Words}: Prime ideals; Cyclic modules; local rings;
 principal ideal rings.}
\thanks {2010{ \it Mathematics Subject Classification}. Primary 13C05, 13E05, 13F10,
  Secondary  13E10, 13H99. }}

\author{{\bf M. Behboodi$^{{\rm a,b}}$\thanks{Corresponding
author.} and} {\bf A. Moradzadeh-Dehkordi$^{{\rm a}}$}\\
 {\small{ $^{\rm a}$Department of Mathematical Sciences, Isfahan University of Technology}}\vspace{-1mm}\\
  {\small{ P.O.Box: 84156-83111, Isfahan, Iran}}\\
 {\small{ $^{\rm b}$School of Mathematics, Institute for Research in Fundamental Sciences
 (IPM)}}\vspace{-1mm}\\ {\small{ P.O.Box: 19395-5746, Tehran, Iran}}\vspace{-1mm}\\
 {\small{mbehbood@cc.iut.ac.ir}}\vspace{-1mm}\\
{\small{a.moradzadeh@math.iut.ac.ir}}}
\date{}
\begin{document}
\maketitle
\begin{abstract}
\noindent {In this paper we study commutative rings $R$ whose
prime ideals are direct sums of cyclic modules. In the case $R$ is
a finite direct product of commutative local rings, the structure
of such rings is completely described. In particular, it is shown
that for a local ring $(R, \cal{M})$, the  following statements
are  equivalent:
  (1)  Every prime ideal of $R$ is a direct sum of cyclic
$R$-modules;
 (2) ${\cal{M}}=\bigoplus_{\lambda\in \Lambda}Rw_{\lambda}$ and $R/{\rm Ann}(w_{\lambda})$
 is a principal ideal ring for each $\lambda \in \Lambda$;
  (3) Every prime ideal of $R$ is a direct sum of at most
$|\Lambda|$ cyclic $R$-modules; and  (4)  Every prime ideal of $R$
is a summand of a direct sum of cyclic $R$-modules. Also, we
establish a theorem which state that, to check whether every prime
ideal in a Noetherian local ring $(R, \cal{M})$ is a direct sum of
 (at most $n$) principal ideals, it suffices to test only the maximal ideal
$\cal{M}$. }
\end{abstract}

{\bf ~~~~~~~~~~~~~~~~~~~~~~~~~~~~~~~~1. Introduction}

It was shown by K\"{o}the \cite{Kothe} that an Artinian
commutative ring $R$ has the property that every module is a
direct sum of cyclic modules if and only if $R$ is  a  principal
ideal ring. Later Cohen-Kaplansky \cite{Cohen2} obtained the
following result: {\it ``a commutative ring $R$ has the property
that every module is a direct sum of cyclic modules if and only if
$R$ is an Artinian principal ideal ring}." (Recently, a
generalization of the K\"{o}the-Cohen-Kaplansky theorem have been
given  by Behboodi et al., \cite{Behboodi} for the noncommutative
setting.) Therefore, an interesting natural question of this sort
is {\it ``Whether the same is true if  one only assumes that every
ideal of $R$ is a direct sum of cyclic modules?"} More recently,
this question was answered by Behboodi et al. \cite{Behboodi1} and
\cite{Behboodi0} for the case $R$ is a finite direct product of
commutative  local rings.

 We note that two theorems from commutative algebra due to I. M. Isaacs and I. S. Cohen state that,
to check whether every ideal in a ring is cyclic (resp. finitely
generated), it suffices to test only the prime ideals (see
\cite[p. 8, Exercise 10]{Kaplansky2} and \cite[Theorem
2]{Cohen1}). So this raises the natural question: "If every prime
ideal of R is a direct sum of cyclics, can we conclude that all
ideals are direct sums of cyclics?" This is not true in general.
In \cite[Example 3.1]{Behboodi1}, for each integer $n\geq 3$, we
provide an example of an Artinian local ring $R$ such that every
prime ideal of $R$ is a direct sum of cyclic $R$-modules, but
there exists a two generated ideal of $R$ which is not a direct
sum of cyclic $R$-modules. Therefore, another  interesting natural
question of this sort is {\it  ``What is the class of commutative
rings $R$ for which every prime ideal is a direct sum of cyclic
modules?"} The goal of this paper is to answer this question in
the case $R$ is a finite direct product of commutative local
rings.  The structure of such rings is completely described.

Throughout this paper, all rings are  commutative  with identity
and all modules are unital. For a ring $R$ we denote by Spec$(R)$
and Max$(R)$
 for the set of prime ideals and maximal ideals  of $R$, respectively. We denote the classical Krull dimension
of $R$ by dim$(R)$. Let $X$ be either an element or a subset of
$R$. The  annihilator of $X$ is the ideal Ann$(X)=\{a\in R~|~
aX=0\}$. A ring $R$  is {\it local} (resp. {\it semilocal}) in
case $R$ has a unique maximal ideal  (resp. a finite number of
maximal ideals). In this paper $(R,{\mathcal{M}})$ will be a local
ring with maximal ideal ${\mathcal{M}}$. A non-zero $R$-module $N$
is called {\it simple} if it has no submodules except $(0)$ and
$N$.

  For a ring $R$, it is shown that if every
prime ideal of $R$ is a direct sum of cyclic $R$-modules, then
${\rm dim}(R)\leq 1$ (Proposition 2.1). Let  $R$ be a semilocal
ring such that every prime ideal of $R$ is a direct sum of cyclic
$R$-modules. Then: (i) $R$ is a principal ideal ring if and only
if every maximal ideal of $R$ is principal (Theorem 2.4); (ii)
 $R$ is a Noetherian ring if and only if every maximal ideal of $R$ is
finitely generated (Theorem 2.5). Also, in Proposition 2.6, it is
shown that  if for each ${\cal{M}}\in {\rm Max(R)}$,
${\cal{M}}=\oplus_{\lambda\in \Lambda}Rw_{\lambda}$ such that for
each $\lambda \in \Lambda$, $R/{\rm Ann}(w_{\lambda})$ is a
principal ideal ring,  then every prime ideal of $R$ is a direct
sum of cyclic modules. However Example 2.7 shows that the converse
 is not true in general, but it is  true when
$R$ is a local ring (see Theorem 2.10). In particular, in Theorem
2.10,  we show that for a local ring $(R, \cal{M})$ the following
statements are equivalent:\vspace{1mm}\\
(1)  Every prime ideal of $R$ is a direct sum of cyclic
$R$-modules.\\
(2) ${\cal{M}}=\bigoplus_{\lambda\in \Lambda}Rw_{\lambda}$ and $R/{\rm Ann}(w_{\lambda})$ is a principal ideal ring for each $\lambda \in \Lambda$.\\
(3) Every prime ideal of $R$ is a
direct sum of at most $|\Lambda|$ cyclic $R$-modules.\\
(4)  Every prime ideal of $R$ is a summand of a direct sum of
cyclic $R$-modules.

Also, if $(R, \cal{M})$  is   Noetherian,   we show that the above
conditions are also equivalent to: (5)  {\it ${\cal{M}}$ is a
direct sum of cyclic $R$-modules} (see Theorem 2.12);  which state
that, to check whether every prime ideal in a Noetherian local
ring $(R, \cal{M})$ is a direct sum of
 (at most $n$) principal ideals, it suffices to test only the maximal ideal
$\cal{M}$.

 Finally, as a consequence, we obtain: if $R=R_1\times
\cdots\times R_k$, where each $R_i$ $(1\leq i\leq k)$ is a local
ring, then  every prime ideal of $R$ is a direct sum of cyclic
$R$-modules if and only if each $R_i$ satisfies the above
equivalent conditions (see Corollary 2.14). We note that the
corresponding result in the case
$R=\prod_{\lambda\in\Lambda}R_\lambda$ where $\Lambda$ is an
infinite index set and each $R_\lambda$ is a local ring, is not
true in general (see Example 2.15). \\

\noindent{\bf 2. Main Results}\vspace{3mm}

We begin with the following evident useful proposition (see
\cite[Proposition 2.5]{Behboodi0}).

\noindent{\bf Proposition 2.1.} {\it Let $R$  be a  ring. If every
prime ideal of $R$ is a direct sum of cyclic $R$-modules, then for
each prime ideal $P$ of $R$, the ring $R/P$ is a principal ideal
domain. Consequently,  $dim(R)\leq 1$.}

\noindent{\bf Proof.} Assume that every prime ideal of $R$ is a
direct sum of cyclic $R$-modules and $P\subseteq Q$ are prime
ideals of $R$. Since  $Q$ is  direct sum of cyclics, we conclude
that $Q/P$ is a principal ideal of $R/P$. Thus every prime ideal
of the ring $R/P$ is  principal and hence  by Lemma 2.1, $R/P$ is
a PID. Since   this holds for all prime ideals $P$ of $R$, thus
dim$(R)\leq 1$.$~\square$

The following two famous theorems from commutative algebra are
crucial in our investigation.

\noindent{\bf Lemma 2.2.} (Cohen \cite[Theorem 2]{Cohen1}) {\it
Let $R$ be a commutative ring. Then $R$ is a Noetherian ring if
and only if every prime ideal of $R$ is finitely generated.}

\noindent{\bf Lemma 2.3.} (Kaplansky \cite[Theorem
12.3]{Kaplansky2})
 {\it A commutative Noetherian ring $R$ is a principal ideal ring if and
only if every maximal ideal of $R$ is principal.}

The following theorem is an analogue of Kaplansky's theorem.

\noindent{\bf Theorem 2.4.} {\it Let $R$  be a semilocal ring such
that every prime ideal of $R$ is a direct sum of cyclic
$R$-modules. Then $R$ is a principal ideal ring if and only if
every maximal ideal of $R$ is principal.}

\noindent{\bf Proof.} $(\Rightarrow)$ is clear.\\
 $(\Leftarrow)$.  We can write $R=R_1\times\ldots\times R_n$ where each $R_i$ is an indecomposable ring
 (i.e., $R_i$ has no any nontrivial idempotent elements).  Clearly  every prime ideal of $R$ is a
 direct sum of cyclic $R$-modules  if and only if   every prime ideal of $R_i$ is a
 direct sum of cyclic $R$-modules for each $1\leq i\leq n$. Also,   every maximal ideal of $R$ is
  principal if and only if
  every maximal ideal of $R_i$ is principal  for each $1\leq i\leq n$.
  Thus without loss of generality, we can assume that $R$ is an indecomposable
  ring. Also, by Proposition 2.1, dim$(R)\leq
1$. \indent  Suppose, contrary to our
 claim, that $R$ is not a principal ideal ring. Thus by Lemma 2.2, $R$ is not a Noetherian ring.
  Thus by Lemma 2.2, there exists a prime ideal $P$  of $R$ such that $P=\bigoplus_{\lambda\in \Lambda}
Rw_{\lambda}$ where $\Lambda$ is an infinite index set and $0\neq
w_{\lambda}\in R$ for each ${\lambda\in \Lambda}$. Thus $P$ is not
a maximal ideal of $R$ and so it is a minimal prime ideal of $R$.\\
\indent For each ${\lambda\in \Lambda}$, there exists a maximal
submodule  $K_{\lambda}$ of $Rw_{\lambda}$ and so
Ann$(Rw_{\lambda}/K_{\lambda})={\cal M}$  for some maximal ideal
${\cal M}$ of $R$. Since ${\rm Max}(R)$ is finite and
$|\Lambda|=\infty$, we can assume that $\{1,2\}\subseteq \Lambda$
and  there exists  ${\cal M}\in {\rm Max}(R)$ such that
$${\rm Ann} (Rw_{1}/K_{1})={\cal M}={\rm Ann} (Rw_{2}/K_{2}).$$
\indent Now set $P=Rw_1\oplus Rw_2\oplus L$ where $L$ is an ideal
of $R$ and $\bar{R}:=R/(K_1\oplus K_2\oplus L)$. Since
${\cal{M}}(Rw_i/K_i)=(0)$
 for $i=1, 2$ and
$$\bar{P}=P/(K_1\oplus K_2\oplus L)\cong (Rw_1\oplus
Rw_2)/(K_1\oplus K_2)\cong R/{\cal{M}}\oplus R/{\cal{M}}$$
 we conclude that ${\bar{{\cal{M}}}}{\bar{P}}=(0)$. It
follows that ${\bar{P}}$ is the only non-maximal prime ideal of
${\bar{R}}$.
 Thus by Lemma 2.2, ${\bar{R}}$ is a Noetherian ring (since ${\bar{P}}$ is  finitely
 generated and every maximal ideal of ${\bar{R}}$ is cyclic) and so by
Lemma 2.3, ${\bar{R}}$ is a principal ideal ring. But ${\bar{P}}$
is a direct sum of two isomorphic simple $R$-modules (so
${\bar{P}}$ is a 2-dimensional  $R/{\cal{M}}$-vector space)
     and hence  it is
not a cyclic $R$-module, a contradiction.$~\square$

Also, the following result is an analogue of Cohen's theorem.

\noindent{\bf Theorem 2.5.} {\it Let $R$  be a semilocal ring such
that every prime ideal of $R$ is a direct sum of cyclic
$R$-modules. Then $R$ is a Noetherian ring if and only if every
maximal ideal of $R$ is finitely generated.}

\noindent{\bf Proof.} $(\Rightarrow)$ is clear.\\
 $(\Leftarrow)$.  We can write $R=R_1\times\ldots\times R_n$ where each $R_i$ is an indecomposable ring
 (i.e., $R_i$ has no any nontrivial idempotent elements).
  Thus without loss of generality, we can assume that $R$ is an indecomposable
  ring with maximal ideals ${\cal{M}}_1, {\cal{M}}_2,\ldots, {\cal{M}}_k$.
  Then by Proposition 2.1, dim$(R)\leq
1$. Suppose, contrary to our
 claim, thus by Lemma 2.2, there exists a prime ideal $P$  of $R$ such that $P=\bigoplus_{\lambda\in \Lambda}
Rw_{\lambda}$ where $\Lambda$ is an infinite index set and $0\neq
w_{\lambda}\in R$ for each ${\lambda\in \Lambda}$. Thus $P$ is not
a maximal ideal of $R$ and so it is a minimal prime ideal of $R$.
Also, by hypothesis for each $1\leq i\leq k$, there exist
$x_{i1},\ldots,x_{in_i}\in R$ such that
$${\cal M}_i=Rx_{i1}\oplus Rx_{i2}\oplus \ldots\oplus
Rx_{in_i}.$$  Since $P$ is a non-maximal prime ideal, without loss
of generality, we can assume that $x_{11},
x_{21},\ldots,x_{k1}\notin P$. It follows that  $Rx_{i2}\oplus
\ldots\oplus Rx_{in_i}\subseteq P$ for each $i=1,\ldots,k$. Set
$$L=(Rx_{12}\oplus  \ldots\oplus Rx_{1n_1})+(Rx_{22}\oplus
\ldots\oplus Rx_{2n_2})+\ldots+(Rx_{k2}\oplus  \ldots\oplus
Rx_{kn_k})$$ Then $L\subseteq P$ and so $L\subseteq
\bigoplus_{\lambda \in {{\Lambda}'}} Rw_{\lambda}$ where $\Lambda
'$ is a finite subset of $\Lambda$.\\
 \indent Clearly, for each ${\lambda\in \Lambda}$,
there exists a maximal submodule  $K_{\lambda}$ of $Rw_{\lambda}$
and hence   Ann$(Rw_{\lambda}/K_{\lambda})={\cal M}$ for some
maximal ideal ${\cal M}$ of $R$.
 Since ${\rm Max}(R)$ is finite and $|\Lambda|=\infty$, we can assume that
$\{1,2\}\subseteq \Lambda$ and  there exists  ${\cal M}\in {\rm
Max}(R)$ such that
$${\rm Ann} (Rw_{1}/K_{1})={\cal M}={\rm Ann} (Rw_{2}/K_{2}).$$
\indent Now we can assume that $P=Rw_1\oplus Rw_2\oplus L$ such
that $\bigoplus_{\lambda \in {{\Lambda}'}} Rw_{\lambda} \subseteq
L$.
 Set $$\bar{R}=R/(K_1\oplus K_1\oplus L).$$ Since
${\cal{M}}(Rw_i/K_i)=(0)$
 for $i=1, 2$ and
$$\bar{P}=P/(K_1\oplus K_2\oplus L)\cong (Rw_1\oplus
Rw_2)/(K_1\oplus K_2)\cong R/{\cal{M}}\oplus R/{\cal{M}},$$ we
conclude that ${\bar{{\cal{M}}}}{\bar{P}}=(0)$. It follows that
${\bar{P}}$ is the only non-maximal prime ideal of ${\bar{R}}$. On
the other hand, for each $1\leq i\leq k$, $Rx_{i2}\oplus
\ldots\oplus Rx_{in_i}\subseteq \bigoplus_{\lambda \in
{{\Lambda}'}} Rw_{\lambda} \subseteq L$. Thus we conclude that
every  maximal ideal of ${\bar {R}}$ is cyclic. Thus by Theorem
2.4, ${\bar{R}}$ is a principal ideal ring.
  But ${\bar{P}}$ is a direct sum of two isomorphic simple $R$-modules (so
   ${\bar{P}}$ is a 2-dimensional  $R/{\cal{M}}$-vector space),
     and hence  it is
not a cyclic $R$-module, a contradiction.$~\square$

\noindent{\bf Proposition 2.6.} {\it Let $R$  be a ring. If for
each ${\cal{M}}\in Max(R)$, ${\cal{M}}=\bigoplus_{\lambda\in
\Lambda}Rw_{\lambda}$ such that for each $\lambda \in \Lambda$,
$R/Ann(w_{\lambda})$ is a principal ideal ring, then every prime
ideal of $R$ is a direct sum of cyclic modules.}

 \noindent{\bf Proof.} Assume that $P$ is a non-maximal prime ideal of $R$. There exists a maximal ideal
 ${\cal{M}}\in {\rm Max(R)}$ such that $P\subsetneqq {\cal{M}}=\bigoplus_{\lambda\in
\Lambda}Rw_{\lambda}$. Thus there exists a $\lambda_0 \in \Lambda$ such
 that $w_{\lambda_0}\notin P$.
 Thus,  $\bigoplus_{\lambda\in {\Lambda}\setminus\{\lambda_0\}}Rw_{\lambda}\subseteq P$ and so by modular
 property, we have
 \begin{center}
  $P=P\cap {\cal{M}}=(P\cap Rw_{\lambda_0})\oplus
  (\bigoplus_{\lambda\in {\Lambda}\setminus\{\lambda_0\}}Rw_{\lambda}).$
  \end{center}
 Now since $P\cap Rw_{\lambda_0}\subseteq Rw_{\lambda_0} \cong R/{\rm Ann}(Rw_{{\lambda}_0})$
  and $R/{\rm Ann}(Rw_{{\lambda}_0})$ is a principal ideal ring, we conclude that  $P\cap Rw_{\lambda_0}$ is
   cyclic. Therefore, $P$  is a direct sum of cyclic modules.~$~\square$

However the following example shows that the converse of
Proposition 2.6, is
 not true in
 general, but we will show  in Theorem 2.10, it is  true when $R$ is a local ring.

\noindent{\bf Example 2.7.}  Let  $R$ be the
 subring  of all sequences from the ring $\prod_{i\in\Bbb{N}}{\Bbb Z}_2$ that are eventually constant.
Then $R$ is a zero-dimensional Boolean ring with minimal prime
ideals $P_i= \{\{a_n\}\in R~|~a_i=0\}$ and
$P_{\infty}=\{\{a_n\}\in R~|~a_n=0~ for~ large~ n\}$ (See
\cite{Anderson}).  Clearly, each $P_i$ is cyclic (in fact
$P_i=Rv_i$ where $v_i=(1,1,\cdots,1,0,1,1,\cdots)$)  and
$P_{\infty}=\bigoplus_{i\in\Bbb{N}}{\Bbb
Z_2}=\bigoplus_{i\in\Bbb{N}}Rw_i$ where
$w_i=(0,0,\cdots,0,1,0,0,\cdots)$. Thus every prime ideal of $R$
is a direct sum of cyclic modules. But the factor ring $R/{\rm
Ann}(v_1)= R/{\rm Ann}(0,1,1,1,\cdots)$ is not a principal ideal
ring (since prime ideal $P_{\infty}/{\rm Ann }(v_1)$ is not a
principal ideal of $R/{\rm Ann}(v_1)$). Also, one can easily to
see that if $P_1=\bigoplus_{\lambda\in\Lambda} Rz_\lambda$ where
$\Lambda$ is an index set and $z_\lambda \in P_1$, then
$|\Lambda|=1$ and $P_1=Rz_\lambda=Rv_1$. Thus the converse of
Proposition 2.6 is not true in general.

By using  Nakayama's lemma, we obtain the following
lemma.

\noindent{\bf Lemma 2.8.} {\it Let $R$ be a ring and $M$ be an
$R$-module such that $M$ is a direct sum of a family of  finitely
generated $R$-modules.  Then Nakayama's lemma holds for $M$ (i.e.,
for each $I\subseteq J(R)$, if $IM=M$, then $M=(0)$).}

 \noindent{\bf Lemma 2.9.} (See Warfield \cite[Proposition 3]{Warfield1}) {\it Let $R$ be a local ring and
 $N$  an  $R$-module. If $N=\oplus_{\lambda\in
\Lambda} R/I_{\lambda}$  where each $I_{\lambda}$ is an ideal of
$R$, then every summand of $N$  is also a direct sum of cyclic
$R$-modules,
 each isomorphic to one of the $R/I_{\lambda}$.}

The following main theorem is an answer to the question  "What is
the class of local rings R for which every prime ideal is a direct
sum of cyclic modules?"

\noindent{\bf Theorem 2.10.} {\it Let $(R, \cal{M})$  be a local
ring. Then the following statements are equivalent:}\vspace{2mm}\\
(1) {\it Every prime ideal of $R$ is a direct sum of cyclic
$R$-modules.}\vspace{1mm}\\
(2) {\it ${\cal{M}}=\bigoplus_{\lambda\in \Lambda}Rw_{\lambda}$
and $R/Ann(w_{\lambda})$
 is a principal ideal ring for each $\lambda \in \Lambda$.}\vspace{1mm}\\
(3) {\it Every prime ideal of $R$ is a
direct sum of at most $|\Lambda|$ cyclic $R$-modules.}\vspace{1mm}\\
(4) {\it Every prime ideal of $R$ is a summand of a direct sum of cyclic
$R$-modules.}

 \noindent{\bf Proof.} $(1) \Rightarrow (2)$. First, we assume that ${\cal{M}}$ is cyclic and  so
 ${\cal{M}}=Rx$ for some   $x\in {\cal{M}}$.  If
 Spec$(R)=\{{\cal{M}}\}$, then by Lemma 2.2,  $R$ is a Noetherian ring and by Lemma 2.3,  $R$ is a
 principal ideal ring. Therefore, $R/ {\rm {Ann}}(x)$ is a principal ideal ring.
 If  Spec$(R)\neq \{{\cal{M}}\}$,
 then for each non-maximal  prime ideal $P$ of $R$,  $x\notin P\subsetneqq {\cal{M}}$. Thus $Px=P$ and so by
  Lemma 2.8, $P=0$. Thus $R$ is a principal ideal domain and so $R/ {\rm {Ann}}(x)$ is principal ideal ring.\\
 \indent Now assume that ${\cal{M}}$ is not cyclic. Then  by hypothesis
  ${\cal{M}}=\bigoplus_{\lambda\in \Lambda}Rw_{\lambda}$
  such that $\Lambda$ is an index set with $|\Lambda|\geq 2$ and $0\neq w_{\lambda}\in {\cal{M}}$ for each $\lambda\in\Lambda$.
   If Spec$(R)=\{{\cal{M}}\}$, then the only maximal ideal of $R/{\rm Ann}(w_{\lambda})$ is principal for each $\lambda\in \Lambda$.
    Thus by Lemma 2.2, $R/{\rm Ann}(w_{\lambda})$ is a Noetherian ring and so by Lemma 2.3, $R/{\rm Ann}(w_{\lambda})$
    is a principal ideal ring for each $\lambda\in \Lambda$.
 If Spec$(R)\neq \{{\cal{M}}\}$, then for each non-maximal  prime ideal $P$  of $R$,  there exists ${\lambda_0}\in \Lambda$ such that  $w_{\lambda_0}\notin
  P$. It follows that  $\bigoplus_{\lambda\in {\Lambda}\setminus\{{\lambda}_0\}}Rw_{\lambda}\subseteq P$.
   Now  by modular property we have
 \begin{center}
  $P=P\cap {\cal{M}}=(P\cap Rw_{\lambda_0})\oplus(\bigoplus_{\lambda\in {\Lambda}\setminus\{{\lambda}_0\}}Rw_{\lambda}).$
  \end{center}
  It follows that $Pw_{\lambda_0}=(P\cap Rw_{\lambda_0})w_{\lambda_0}$. Also  since $\lambda_0\notin P$,
$P\cap Rw_{\lambda_0}=Pw_{\lambda_0}$ and hence
$Pw_{\lambda_0}=(Pw_{\lambda_0})Rw_{\lambda_0}$. Now by Lemma 2.8,
  $Pw_{\lambda_0}=0$, since $Pw_{\lambda_0}$ is a direct sum of cyclic $R$-modules.
   Therefore, $P=\bigoplus_{\lambda\in {\Lambda}\setminus\{{\lambda}_0\}}Rw_{\lambda}$.
   Thus we conclude that ${\rm Spec}(R)=\{ {\bigoplus_{{\lambda\in {{\Lambda}\setminus\{\lambda}_j\}} }{Rw_{\lambda}}~|~{w_{{\lambda}_j}}}\notin {\rm Nil}(R) \}$.
   This shows that  for each $\lambda\in \Lambda$, all  prime ideals
  of $R/{\rm Ann}(w_{\lambda})$ are principal. Thus by Lemmas 2.2 and 2.3, $R/{\rm Ann}(w_{\lambda})$ is a principal ideal ring for each $\lambda\in \Lambda$.\\
 $(2) \Rightarrow (3)$. Assume that ${\cal{M}}$ is  cyclic and  so ${\cal{M}}=Rx$ for some $x\in {\cal{M}}$. If
 Spec$(R)=\{{\cal{M}}\}$, then the proof is complete. If  Spec$(R)\neq \{{\cal{M}}\}$,
 then for each non-maximal  prime ideal $P$ of $R$,  $x\notin P\subsetneqq{\cal{M}}$. Thus $Px=P$ and so $Px=Px(Rx)$.
 By hypothesis $R/ {\rm {Ann}}(x)$ is a principal ideal ring and so $Px$ is principal. Thus by Lemma 2.8, $Px=0$ and so $P=0$.\\
   \indent Now assume that ${\cal{M}}$ is not cyclic and so ${\cal{M}}=\bigoplus_{\lambda\in \Lambda}Rw_{\lambda}$ such that
    $|\Lambda|\geq 2$ and $R/{\rm Ann}(Rw_{\lambda})$ is a principal ideal ring for each $\lambda \in \Lambda$. If Spec$(R)=\{{\cal{M}}\}$,
    then the proof is complete.
   If Spec$(R)\neq \{{\cal{M}}\}$, then for each non-maximal prime ideal $P$  of $R$,  there exists
   ${\lambda_0}\in \Lambda$ such that  $w_{\lambda_0}\notin
  P$. This implies that $\bigoplus_{\lambda\in {\Lambda}\setminus\{{\lambda}_0\}}Rw_{\lambda}\subseteq P$.
   Thus by modular property,
  $P=P\cap {\cal{M}}=(P\cap Rw_{\lambda_0})\oplus (\bigoplus_{\lambda\in {\Lambda}\setminus\{{\lambda}_0\}}Rw_{\lambda})$
  and so $Pw_{\lambda_0}=(P\cap Rw_{\lambda_0})w_{\lambda_0}$. Also $P\cap Rw_{\lambda_0}=Pw_{\lambda_0}$, since $\lambda_0\notin P$ and  $Pw_{\lambda_0}=Pw_{\lambda_0}(Rw_{\lambda_0})$. But $Pw_{\lambda_0}=P\cap Rw_{\lambda_0}$ is principal, since
  $R/{\rm Ann}(w_{{\lambda}_0})$ is a principal ideal ring. Thus by Lemma 2.8, $Pw_{\lambda_0}=0$ and
  so  $P=\bigoplus_{\lambda\in {\Lambda}\setminus\{{\lambda}_0\}}Rw_{\lambda}$. Thus we conclude that
     \begin{center}
    ${\rm Spec}(R)=\{{\cal M}\}\cup\{ {\bigoplus_{{\lambda\in {{\Lambda}\setminus\{\lambda}_j\}} }{Rw_{\lambda}}~|~{w_{{\lambda}_j}}}\notin {\rm Nil}(R) \}$,
    \end{center}
    and hence every prime ideal of $R$ is a
direct sum of at most $|\Lambda|$ cyclic $R$-modules.\\
       $(3) \Rightarrow (4)$ is clear.\\
   $(4) \Rightarrow (1)$ is  by Lemma 2.9.~$~\square$

Also, the following result  is an answer to the question "What is
the class of local rings $(R, \cal{M})$ for which $\cal{M}$ is
finitly generated and every prime ideal is a direct sum of cyclic
modules?"

\noindent{\bf Corollary 2.11} {\it Let $(R, \cal{M})$  be a local
ring. Then the following statements are equivalent:}\vspace{2mm}\\
(1) {\it $R$ is a Noetherian ring and every prime ideal of $R$ is a direct sum of cyclic
$R$-modules.}\vspace{1mm}\\
(2) {\it ${\cal{M}}=\bigoplus_{i=1}^{n}Rw_{i}$ and $R/Ann(w_{i})$ is a principal ideal ring for each  $1\leq i\leq n$.}\vspace{1mm}\\
(3) {\it Every prime ideal of $R$ is a
direct sum of at most $n$ cyclic $R$-modules.}\vspace{1mm}\\
(4) {\it $R$ is a Noetherian ring and every prime ideal of $R$ is a summand of a direct sum of \indent cyclic
$R$-modules.}

 \noindent{\bf Proof.} The proof is straightforward by   Theorem 2.5 and Theorem 2.10.$~\square$

Next, we greatly improve the main theorem above (Theorem 2.10) in
the case $R$
 is a Noetherian local ring. In fact, we establish the
following result which state that, to check whether every prime
ideal in a Noetherian local ring $(R, \cal{M})$ is a direct sum of
 (at most $n$) principal ideals, it suffices to test only the maximal ideal
$\cal{M}$. We note that this is also  a generalization of the
Kaplansky Theorem in the case $R$ is a Noetherian local ring.

\noindent{\bf Theorem  2.12.} {\it Let $(R, \cal{M})$  be a
Noetherian local
ring. Then the following statements are equivalent:}\vspace{2mm}\\
(1) {\it Every prime ideal of $R$ is a direct sum of cyclic
$R$-modules.}\vspace{1mm}\\
(2) {\it ${\cal{M}}=\bigoplus_{i=1}^{n}Rw_{i}$ and $R/Ann(w_{i})$ is a principal ideal ring for each  $1\leq i\leq n$.}\vspace{1mm}\\
(3) {\it The maximal ideal ${\cal{M}}$ is a direct sum of $n$
cyclic $R$-modules.}\vspace{1mm}\\
(4) {\it Every prime ideal of $R$ is a direct sum of at most $n$ cyclic $R$-modules.}\vspace{1mm}\\
(5) {\it Every prime ideal of $R$ is a summand of a direct sum of cyclic $R$-modules.}

 \noindent{\bf Proof.} $(1) \Rightarrow (2)$ and  $(4) \Rightarrow (5) \Rightarrow (1)$
  are by Theorem 2.10.\\
   $(2) \Rightarrow (3)$  is clear. \\
$(3) \Rightarrow (4)$. If ${\cal{M}}=Rx$ is a cyclic $R$-module,
then by Lemma 2.3, $R$ is a principal ideal ring. Assume that
${\cal{M}}=\bigoplus_{i=1}^n Rw_i$ where $n\geq 2$. If
Spec$(R)=\{\cal{M}\}$, then the proof is compleat. Thus we can
assume that  Spec$(R)\neq \{\cal{M}\}$ and suppose that
$P\subsetneqq\cal{M}$ is a prime ideal of $R$. Without loss of
generality, we can assume that, $w_1\notin P$. This implies that
$\bigoplus_{i=2}^nRw_{i}\subseteq P$. Now  by modular property we
have  $P=P\cap {\cal{M}}=(P\cap Rw_1)\oplus
(\bigoplus_{i=2}^nRw_{i})$,
 and hence  $Pw_1=(P\cap
Rw_1)w_1$. Also,  $P\cap Rw_1=Pw_1$ since $w_1\notin P$. Thus
$Pw_1=(Pw_1)Rw_1$,
 and so  by Lemma 2.8, $Pw_1=0$. Therefore,
$P=\bigoplus_{i=2}^nRw_{i}$.~$~\square$

 \noindent{\bf Remark  2.13.}  Let   $R=R_1\times \cdots\times
R_k$ where $k\in\Bbb{N}$ and each $R_i$ is a nonzero ring. One can
easily see that, each prime ideal $P$  of $R$ is of the form
$P=R_1\times \cdots\times R_{i-1}\times P_i\times
R_{i+1}\times\cdots\times R_k$ where $P_i$ is a prime ideal of
$R_i$. Also, if $P_i$ is a direct sum of $\Lambda$ principal
ideals of $R_i$, then it is easy to see that $P$ is also a direct
sum of $\Lambda$ principal ideals of $R$. Thus the ring $R$ has
the property that whose prime ideals are direct sum of cyclic
$R$-modules if and only if for each $i$ the ring $R_i$ has this
property.

We are thus led to the following
strengthening of Theorem 2.10.

\noindent{\bf Corollary  2.14.} {\it Let $R=R_1\times \cdots\times
R_k$ where $k\in\Bbb{N}$ and each $R_i$ is  a  local ring with
maximal ideal ${\cal{M}}_i$ $(1\leq i\leq k)$. Then the following
statements are equivalent:}\vspace{2mm}\\
(1) {\it Every prime ideal of $R$ is a direct sum of cyclic
$R$-modules.}\vspace{1mm}\\
(2) {\it For each $i$,  ${\cal{M}}_i=\bigoplus_{\lambda_i\in
\Lambda_i}Rw_{{\lambda}_i}$ and
 $R/Ann(w_{{\lambda}_i})$ is a principal ideal ring for each \indent ${\lambda_i} \in {\Lambda_i}$.}\vspace{1mm}\\
(3) {\it Every prime ideal of $R$ is a
direct sum of at most $|\Lambda|$ cyclic $R$-modules, where $\Lambda=\indent max\{\Lambda_i~|~i=1,\ldots,k\}$.}\vspace{1mm}\\
(4) {\it Every prime ideal of $R$ is a summand of a direct sum of
cyclic $R$-modules.}

 \noindent{\bf Proof.} The proof is straightforward by  Theorem 2.10 and Remark
 2.13.$~\square$

We conclude this paper with the following interesting example. In
fact, the following example shows that the corresponding of the
above result in the case $R=\prod_{\lambda\in\Lambda}R_\lambda$
where $\Lambda$ is an infinite index set and each $R_\lambda$ is a
local ring (even if for each $\lambda\in\Lambda$,
$R_\lambda\cong\Bbb{Z}_2$), is not true in general.

  \noindent{\bf Example 2.15.}  Let $R=\prod_{\lambda\in\Lambda}F_\lambda$ be a direct product of fields
  $\{F_\lambda\}_{\lambda\in\Lambda}$ where $\Lambda$ is an infinite index set. Clearly,
   $I=\bigoplus_{\lambda\in\Lambda}F_\lambda$  is a non-maximal  ideal of $R$. Thus there exists a
   maximal ideal $P$ of $R$ such that $I\subsetneqq P$.  It was shown by Cohen and Kaplansky
    \cite[Lemma 1]{Cohen2} that $P$ is not a direct sum of
   principal ideals.

\end{document}